\documentclass[11pt]{article}
\usepackage{latexsym, amsfonts}
\newcommand{\be}{\begin{equation}}
\newcommand{\ee}{\end{equation}}
\newcommand{\bea}{\begin{eqnarray}}
\newcommand{\eea}{\end{eqnarray}}
\newcommand{\bean}{\begin{eqnarray*}}
\newcommand{\eean}{\end{eqnarray*}}
\newcommand{\brray}{\begin{array}}
\newcommand{\erray}{\end{array}}
\newcommand{\ben}{\begin{equation}{nonumber}}
\newcommand{\een}{\end{equation}{nonumber}}

\newtheorem{dfn}{Definition}[section]
\newtheorem{thm}[dfn]{Theorem}
\newtheorem{lmma}[dfn]{Lemma}
\newtheorem{ppsn}[dfn]{Proposition}
\newtheorem{crlre}[dfn]{Corollary}
\newtheorem{xmpl}[dfn]{Example}
\newtheorem{rmrk}[dfn]{Remark}
\newcommand{\bdfn}{\begin{dfn}}
\newcommand{\bthm}{\begin{thm}}
\newcommand{\blmma}{\begin{lmma}}
\newcommand{\bppsn}{\begin{ppsn}}
\newcommand{\bcrlre}{\begin{crlre}}
\newcommand{\bxmpl}{\begin{xmpl}}
\newcommand{\brmrk}{\begin{rmrk}}
\newcommand{\edfn}{\end{dfn}}
\newcommand{\ethm}{\end{thm}}
\newcommand{\elmma}{\end{lmma}}
\newcommand{\eppsn}{\end{ppsn}}
\newcommand{\ecrlre}{\end{crlre}}
\newcommand{\exmpl}{\end{xmpl}}
\newcommand{\ermrk}{\end{rmrk}}
\newcommand{\IC}{\mathbb{C}}
\newcommand{\al}{\alpha}

\newcommand{\gma}{\gamma}

\newcommand{\Gma}{\Gamma}
\newcommand{\cla}{{\cal A}}
\newcommand{\clb}{{\cal B}}

\newcommand{\clh}{{\cal H}}

\newcommand{\clk}{{\cal K}}
\newcommand{\cll}{{\cal L}}

\newcommand{\cls}{{\cal S}}

\newcommand{\clv}{{\cal V}}

\def\a*{{\cal A}_{h,*}}
\def\B{{\cal B}(h)}
\def\B1{{\cal B}_1(h)}
\def\b{{\cal B}^{s. a. }(h)}
\def\b1{{\cal B}^{s. a. }_1(h)}

\newcommand{\ot}{\otimes}

\newcommand{\raro}{\rightarrow}

\def \qed {$\Box$}

\begin{document}
\begin{center}
{\large {\bf On Equivariant Embedding  of  Hilbert $C^*$ modules  }}\\
{\bf Debashish Goswami}\\
\end{center}
AMS subject classification no. : {\bf 46L08}
\begin{abstract}
We prove that an arbitrary (not necessarily countably generated) Hilbert $G$-$\cla$ module on a $G-C^*$ algebra $\cla$ admits an equivariant embedding into a trivial $G-\cla$ module, provided $G$ is a compact Lie group and its action on $\cla$ is ergodic.
\end{abstract}
\section{Introduction}
 Let $G$ be a locally
compact group, $\cla$ be a $C^*$-algebra, and assume that there is a strongly
continuous representation $\al : G \raro Aut(\cla)$.
Following the terminology of \cite{MP}, we introduce the concept of a Hilbert
$C^*$ $G-\cla$-module as follows :\\
 \bdfn {\rm
\label{gamod}
A {\it Hilbert $C^*$ $G-\cla$ module} (or {\it $G-\cla$ module} for short) is a pair
$(E,\beta)$ where $E$ is a Hilbert $C^*$ $\cla$-module and $\beta$ is a map
from $G$ into the set of $\IC$-linear (caution :  {\bf not} $\cla$-linear !)
maps from $E$ to $E$, such that $\beta_g \equiv \beta(g), g \in G$ satisfies
the following :\\
(i) $\beta_{gh}=\beta_g \circ \beta_h$ for $g,h \in G,$ $\beta_e={\rm Id},$ where
$e$ is the identity element of $G$;\\
  (ii) $\beta_g(\xi a)=\beta_g(\xi)
\al_g(a)$ for $\xi \in E, a \in \cla;$\\
(iii) $g \mapsto \beta_g(\xi)$ is
continuous for each fixed $\xi \in E;$\\
(iv) $   \langle   \beta_g(\xi),\beta_g(\eta)  \rangle  =\al_g(  \langle  \xi,\eta  \rangle  )$ for all
$\xi, \eta \in E,$ where $  \langle  \cdot,\cdot  \rangle  $ denotes the $\cla$-valued
inner product of $E.$ } \edfn When $\beta$ is understood from the
context, we may refer to $E$ as a $G-\cla$ module, without
explicitly mentioning the pair $(E,\beta).$ Given two $G-\cla$
 modules $(E_1,\beta)$ and $(E_2, \gamma)$, there is a natural $G$-action
induced on $\cll(E_1,E_2)$, given by
$\pi_g(T)(\xi):=\gamma_g(T(\beta_{g^{-1}}(\xi)))$ for $g \in G, \xi \in E_1, T
\in \cll(E_1,E_2).$ $T \in \cll(E_1,E_2)$ is said to be {\it $G$-equivariant}
if $\pi_g(T)=T {\rm ~for ~all~} g \in G.$ It is clear that for each fixed $T \in
\cll(E_1,E_2)$ and $\xi \in E_1$, $g \mapsto \pi_g(T)\xi$ is continuous. We say
that $T$ is {\it $G$-continuous} if $g \mapsto \pi_g(T) $ is continuous with
respect to the norm topology on $\cll(E_1,E_2).$
 We say that $E_1$ and $E_2$ are {\it isomorphic as $G-\cla$-modules}, or that
they are {\it equivariantly isomorphic} if there is a $G$-equivariant unitary
map $T \in \cll(E_1,E_2).$ We call a $(G-\cla)$ module of the form $(\cla \ot \clh, \al_g \ot \gma_g)$ (where $\clh$ is some Hilbert space) a trivial $G-\cla$ module. We say that $(E,\beta)$ is {\it embeddable } if there is an equivariant isometry from $E$ to $\cla \ot \clh$ for some Hilbert space $\clh$ with a $G$-action $\gma$, or in other words, $(E,\beta)$ is equivariantly isomorphic with a sub-$G-\cla$ module of $(\cla \ot \clh, \beta \ot \gma)$.  Note that $\cla \ot \clh$ is the closure of $\cla \ot_{\rm alg} \clh$  under the norm inherited from $\clb(\clh_0, \clh_0 \ot \clh)$ where $\clh_0$ is any Hilbert space such that $\cla$ is isometrically embedded into $\clb(\clh_0)$. The  following result on the embeddabiity is due to Mingo and Phillips (\cite{MP}).

 \bthm
\label{eqkas}
 Let $(E,\beta)$ be a Hilbert $C^*$
$G-\cla$ module and assume that $E$ is countably generated as a Hilbert
$\cla$-module,   that is, there is a countable set $S=\{e_1,e_2,...\}$ of elements
of $E$ such that the right $\cla$-linear span of $S$ is dense in $E$.
Assume furthermore that $G$ is compact. Then $(E,\beta)$ is embeddable.
  \ethm

  When $G$ is the trivial singleton group, the above result was proved by Kasparov.

   If the $C^*$ algebra $\cla$  is replaced by a  von Neumann algebra $\clb \subseteq \clb(h)$ for
some Hilbert space $h$ and $G$ is a locally compact group with a strongly
continuous unitary representation $g \mapsto u_g \in \clb(h)$, one can define Hilbert von Neumann $G$-$\clb$ module $(E,\beta)$. The only difference is  that $E$ is now a Hilbert von Neumann
$\cla$-module equipped with the natural locally convex strong operator topology, and that we replace the norm-continuity   in (iii) of the above
definition by a weaker continuity : namely, the continuity of $g \mapsto
\beta_g(\xi)$ (for fixed $\xi \in E$) with respect to the locally convex
topology of $E$. In this case, we have  a stronger version of    the Theorem
\ref{eqkas} (see \cite{CGS} and \cite{book}, Theorem 4.3.5, page 99), namely without the condition of $E$ being countably generated and without the compactness of $G$. It should be remarked here that the trivial Hilbert von Neumann $\clb$ module $\clb \ot \clh$ is defined to be the closure of $\clb \ot_{\rm alg} \clh$ with respect to the strong-operator-topology inherited from $\clb(h,h \ot \clh)$.

In Theorem \ref{eqkas}, the assumption that $E$ is countably generated  restricts the applicability of the result, since it is not always easy to check the property of being countably generated. However, under some special assumption on the $G$-$C^*$ algebra $\cla$, i.e. conditions on the group $G$, the $C^*$ algebra  $\cla$ and also on the nature of the action, it may be possible to prove the embedability for an arbitrary Hilbert $G$-$\cla$ module. The aim of the present article is to give some such sufficient conditions.
\section{Ergodic action and its implication}
We say that the action $\al$ of $G$ on a unital $C^*$-algebra
$\cla$ is {\it ergodic} if $\al_g(a)=a ~{\rm ~for ~all~} g \in G$ if and
only if $a$ is a scalar multiplie of $1$. There is a considerable amount of literature  on  ergodic action of compact groups, and we shall
quote one interesting structure theorem  which
will be useful for  us. \\
\bppsn
\label{ergaction}
%By combining the results of
%\cite{Shiga}, \cite{Alb} and \cite{HKr}
Let $G$ be a compact group acting ergodically on a unital $C^*$-algebra
$\cla.$ Then there is a set of elements  $t^\pi_{ij},\pi \in
\hat{G},i=1,...,d_\pi,j=1,...,m_\pi$ of  $\cla$, where $\hat{G}$ is the set of
 equivalence classes of irreducible representations  of $G$, $d_\pi$ is the dimension of the
irreducible representation  space denoted by $\pi$, $m_\pi \leq d_\pi$ is a
natural number,  such that the followings hold : \\
(i) There is a unique faithful $G$-invariant state $\tau$ on $\cla$, which is in fact
a trace,\\
  (ii) The linear span of $\{ t^\pi_{ij} \}$ is norm-dense in
 $\cla$,\\
 (iii) $\{ t^\pi_{ij} \}$ is an orthonormal basis of
 $h=L^2(\cla,\tau)$, \\
 (iv) The action of $u_g$ coincides with the $\pi$th irreducible
 representation of $G$ on the vector space spanned by
 $t^\pi_{ij},i=1,...,d_\pi$ for each fixed $j$ and $\pi$,\\
 (v) $\sum_{i=1,...d_\pi} (t^\pi_{ij})^*t^\pi_{ik}=\delta_{jk}
 d_\pi 1$, where $\delta_{jk}$ denotes the Kronecker delta symbol.
Thus, in particular, $\| t^\pi_{ij} \| \leq
\sqrt{d_\pi} ~{\rm ~for ~all~} \pi,i,j$. \eppsn
The proof can be obtained
by combining the results of \cite{Shiga},\cite{HKr} and
\cite{Alb}.

Let now  $h=L^2(\cla,\tau)$, where $\tau$ is the unique
invariant faithful trace described in Proposition \ref{ergaction}.
 Let $u_g$ be the unitary in $h$ induced by
the action of $G$,   that  is, on the dense set $\cla \subseteq h,$
$u_g(a):=\al_g(a),$ where $\al_g$ denotes the $G$-action on $\cla$.
Denote also by $\al$ the action $g \mapsto u_g \cdot u_g^*$ on
$\tilde{\cla}$, which is the weak closure of $\cla$ in $\clb(h)$.

Let us now specialize to the case of a compact Lie group. If $G$ is
such a group, with a basis of the Lie algebra given by $\{
\chi_1,...,\chi_N \}$, which has a strongly continuous action
$\theta$ on a Banach  space $F$, we can consider the space of
`smooth' or $C^\infty$-elements of $F$, denoted by $F^\infty$,
consisting of all $\xi \in f$ such that $G \ni g \mapsto
\theta_g(\eta)$ is $C^\infty$. It is easy to prove that  $F^\infty$
is dense in $F$, and it is a $\ast$-subalgebra if $F$ is a locally
convex $\ast$-algebra. Moreover, we equip $E^\infty$ with a family
of seminorms $\| \cdot \|_{\infty,n}$, $n=0,1,...$ given by
$$ \| \xi \|_{\infty,n}:=\sum_{i_1,i_2,...i_k; k \leq n, i_t \in \{
1,...,N \} } \| \partial_{i_1} \partial_{i_2} ... \partial_{i_k} \xi
\|,$$
 with the convention $\| \cdot \|_{\infty,0}=\| \cdot \|$ and where $\partial_j(\xi):=\frac{d}{dt}|_{t=0} \theta_{{\rm exp}(t \chi_j)}(\xi)$. The space $F^\infty$ is complete under this family of seminorms, and thus is a
 Fr$\acute{\rm e}$chet space. When $F$ is  Hilbert space or a Hilbert module, we shall also consider a map $d_j$ given by essentially the same expression as that of $\partial_j$, with $\chi_j$ replaced by $i \chi_j$, and the Hilbertian seminorms $\{ \| \cdot \|_{2,n} \}$ are given by
 $$\| \xi \|^2_{2,n}:=\sum_{i_1,i_2,...i_k; k \leq n, i_t \in \{ 1,...,N \} } \|d_{i_1} d_{i_2} ... d_{i_k} \xi \|_2^2, $$ with $\| \cdot \|_2$ denoting the norm of the Hilbert space (or Hilbert module) $F$.

 More generally, if $F$ is a complete locally convex  space given by a family of seminorms $\{ \| \cdot \|^{(q)} \}$, then we can consider the smooth subspace $F^\infty$ and the maps $\partial_j$ as above, and make it a complete locally convex  space with respect to a larger  family of seminorms $\{ \| \cdot \|^{(q)}_n \}$ where $$ \| \xi \|^{(q)}_n:=\sum_{i_1,i_2,...i_k; k \leq n, i_t \in \{ 1,...,N \} } \| \partial_{i_1} \partial_{i_2} ... \partial_{i_k} \xi \|^{(q)}.$$ In case $F$ is a von Neumann algebra equipped with the locally convex strong operator topology, the locally convex space $F^\infty$ is a topological $\ast$-algebra, strongly dense in $F$.

\blmma \cite{book} \label{automatic_a3} Let $G$ be a compact Lie
group acting ergodically on a unital $C^*$-algebra $\cla.$ Then
$h^{\infty}=\cla^\infty$ as Fr$\acute{\rm e}$chet spaces. \elmma
{\it Proof :}\\
The fact that $\cla^\infty=h^\infty$ as sets is contained in Lemma 8.1.20 of \cite{book} (page 200-201). We only prove that the identity map is a topological homeomorphism.

                   Since
 the trace $\tau$ is finite,  the Fr$\acute{\rm e}$chet topology of $\cla ^\infty$ is
 stronger than that of $h_{\infty}$. This implies that the identity map $I$, viewed as a linear map from the Fr$\acute{\rm e}$chet space $h^\infty$ to the
 Fr$\acute{\rm e}$chet space $\cla^\infty$ is closable, hence continuous. This completes the proof that the two
 Fr$\acute{\rm e}$chet topologies on $\cla^\infty=h^\infty$ are equivalent, i.e. $\cla^\infty=h^\infty$ as topological spaces.   \qed

\blmma \label{mainlemma} Let $\clh$ be a (not necessarily separable)
Hilbert space with a unitary representation $w \equiv w_g$ of $G$,
and let us consider the Fr$\acute{\rm e}$chet modules $(\tilde{\cla}
\ot \clh)^\infty$ and $(\cla \ot \clh)^\infty$ corresponding to the
action $\gamma_g:=\al_g \ot w_g$. Let $\xi$ be an element of
  $(\tilde{\cla} \ot \clh)^\infty$ such that $g \mapsto \gamma_g(\xi)$ is continuous in the operator-norm topology. Then  $\xi$ actually belongs to $ \cla \ot \clh$.
\elmma
{\it Proof :-}\\
We shall denote by $\| \cdot \|_p$ ($p \geq 1$) the $L^p$-norm
coming from the trace $\tau$ on $\cla$. The identity $1$ of $\cla$
will also be viewed as a unit vector in $L^2(\tau)$. Fix an
orthonormal basis $\{ e_\al, \al \in T \}$ of $\clh$ (which need not
be separable), with each $e_\alpha \in \clh^\infty$. Fix $\xi \in
(\tilde{\cla} \ot \clh)^\infty$ satisfying the hytothesis of the
lemma. Since $L^2(\tau)$ is separable, say with an orthonormal basis
given by $\{ x_1,x_2,... \}$, we can find, for each $i$, a counteble
subset $T_i$ of $T$ such that $<\xi 1, x_i \ot e_\alpha>=0$ for all
$\alpha \notin T_i$. Denoting by $T_\infty$ the countable set
$\bigcup_i T_i$, we have   $<\xi1,v \ot e_\alpha>=0~\forall v \in
L^2(\tau)$, for all $\alpha \notin T_\infty$. Write $T_\infty=\{
e_{\alpha_1},e_{\alpha_2},... \}$.  Denote by $\xi_n$ the element in
$\cla^\infty \ot_{\rm alg} \clh^\infty$ given by $\xi_n=(I \ot
P_n)\xi$, where $P_n$ denotes the orthogonal projection  onto the
linear span of $\{ e_{\alpha_1},...,e_{\alpha_n} \}$. It is clear
that $\xi_n1 \raro \xi 1$ as $n \raro \infty$. Now, for a $C^\infty$
complex-valued function $f$ on $G$ and an element $\eta \in
\tilde{\cla} \ot \clh$, denote by $\gamma(f)(\eta)$ the element
$\int_G f(g) \gamma_g(\eta) dg~\in \tilde{\cla} \ot \clh \subseteq
L^2(\tau) \ot \clh$, where $dg$ stands for the normalised Haar
measure on $G$ and the integral is convergent in the strong-operator
topology. We claim that it is enough to prove that $\gamma(f)(\xi)
\in \cla \ot \clh$ for all $f \in C^\infty(G)$. Let us first prove
this claim. Since $g \mapsto \gamma_g(\xi)$ is norm-continuous,
given $\epsilon>0$, we can find a nonempty open subset $U$ of $G$
such that $\| \gamma_g(\xi)-\xi \| \leq \epsilon$ for all $g \in G$,
and then choose $f \in C^\infty(G)$ with ${\rm supp}(f) \subseteq
U$, $f \geq 0$ and $\int_G f dg=1$. It is easy to see that $\|
\gamma(f)(\xi)-\xi \| \leq \epsilon$. Thus, $\xi$ is the
operator-norm limit of a sequence of elements of the form
$\gamma(f)(\xi)$, which proves the claim.

Let us now complete the proof of the lemma by showing that
$\eta:=\gamma(f)(\xi)$ indeed belongs to $\cla \ot \clh$ for every
$f \in C^\infty(G)$. To this end, first observe that
$\eta_n:=\gamma(f)(\xi_n)$ belongs to $\cla^\infty \ot_{\rm alg}
\clh^\infty \subseteq \cla \ot \clh$ for all $n$.  Moreover, since
$u_g1=1$ for all $g$ and $\gamma_g(\cdot)={\rm ad}_{u_g} \ot w_g$,
it is clear that $\eta_n1 \raro \gamma(f)(\xi)1 = \eta 1$ as $n
\raro \infty$.   Since each $\eta_{m,n}$ belongs to $\cla \ot \clh$,
for proving $\eta \in \cla \ot \clh$ it is enough to prove that
$\eta_{m,n} \raro 0$ in the topology of $\cla \ot \clh$, i.e.
$x_{mn}:=<\eta_{m,n}, \eta_{m,n}> \raro 0$ in the norm-topology of
$\cla$.   We shall prove that $x_{mn} \raro 0$ in the Fr$\acute{\rm
e}$chet topology of $h^\infty$, which will prove that it converges
to $0$ also in the topology of $\cla^\infty.$

To this end,
 first note that for  $\beta_1,\beta_2 \in (\tilde{\cla} \ot \clh)$, we have
 $$ \| <\beta_1,\beta_2> \|^2_2 =\tau(\beta_2^* \beta_1 \beta_1^* \beta_2) \leq \| \beta_1 \|^2 \tau(\beta_2^*\beta_2),$$
 hence $\| <\beta_1,\beta_2> \|_2 \leq \| \beta_1 \| \|\beta_2 \|_2$. Moreover,
 $\| <\beta_2,\beta_1> \|_2=\| <\beta_1,\beta_2 >^* \|_2=\|<\beta_1,\beta_2>\|_2$
 (since $\tau(x^*x)=\tau(xx^*),$ we have $\| x\|=\|x^*\|$). From this, we have
 $$ \| < \gamma(f)(\beta), \gamma(f)(\beta)> \|_2 \leq C(f)^2 \| \beta \| \| \beta \|_2,$$ where
 $C(f):=\int |f| dg$. Let us now fix an ordered $k$-tuple $I=(i_1,...,i_k)$ ($k$ nonnegative integer), and
 let $C$ denote the maximum of $C(\chi_{j_1}...\chi_{j_p}f)$ where $J=(j_1,...,j_p)$ varies over
 all (including the empty set) ordered subsets of $I$. Let us abbreviate $\partial_{j_1}...\partial_{j_p}\beta$ and
 $\chi_{j_1}...\chi_{j_p}f$ by $\partial_J\beta$ and $f_J$ respectively,
 for $\beta \in (\cla \ot \clh)^\infty$. Note that $$ \partial_J \gamma(f)(\beta)=(-1)^k\gamma(f_J)(\beta).$$ Using this as well as the Leibniz formula $\partial_I<\beta, \beta>=\sum_J <\partial_J \beta, \partial_{I-J}\beta>$ (with $J$ varying over all ordered substes of $I$), we have the following :
  \bean
 \lefteqn{\| \partial_I x_{mn} \|_2}\\
 &\leq & \sum_J \| <\partial_J (\eta_m-\eta_n), \partial_{I-J}(\eta_m-\eta_n)>\|_2\\
 &\leq & 2^k C^2 \|\xi_m-\xi_n \| \| \xi_m-\xi_n \|_2 \\
 & \leq & 2^{k+1} C^2 \| \xi \| \| \xi_m-\xi_n \|_2, \eean
 since the number of ordered subsets of $I$ is $2^k$ and  it is clear  from the definition of $\xi_n$  that $\| \xi_n \| \leq \|\xi \|$ for all $n$.   We also have $\| \xi_m-\xi_n \|^2_2=\tau(<\xi_m-\xi_n, \xi_m-\xi_n>)=<(\xi_m-\xi_n)(1), (\xi_m-\xi_n)(1)> \raro 0$ as $m,n \raro \infty$. This proves $x_{m,n} \raro 0$ in the topology of $h_\infty$, thereby completing the proof of the lemma.
 \qed

\section{Main results on equivariant embedding of Hilbert modules}
Let $(E,\beta)$ be a $G-\cla$ module, where $\cla$ and $G$ are as in the previous section, i.e. $G$ is a compact Lie group acting ergodically on the $C^*$ algebra $\cla$. In this final section, we shall prove that any such $(E,\beta)$ is embeddable.
\blmma
\label{eqvnkas}
We can find a Hilbert space $\clk$, a strongly continuous unitary representation
 $g \mapsto V_g \in \clb(\clk)$ and a $\cla$-linear isometry $\Gamma_0 : E
\raro \clb(h,\clk)$,  such that $\Gma_0 \beta_g(\xi)=V_g(\Gma_0 \xi)u_g^{-1}$, and moreover,  the complex linear span  of elements of the form $\Gma \xi w$ where $\xi \in E$ and $ w \in h$ is dense in $\clk$.
\elmma
{\it Proof :}\\
The proof of this result is adapted from \cite{CGS} and \cite{book}, Theorem 4.3.5 (page 99-101). We shall give only a brief sketch of the arguments involved, omitting the details. We consider first the formal vector space (say $\clv$) spanned by symbols $(\xi,w)$, with $\xi \in E$ and $w \in h$, and define a semi-inner product on this formal vector space by setting $$ <(\xi,w), (\xi^\prime,w^\prime)>=<w,<\xi,\xi^\prime> w^\prime>,$$ where $<\xi,\xi^\prime>$ denotes the $\cla$-valued inner product on $E$.  By extending this semi-inner product by linearity and then taking quotient by the subspace (say $\clv_0$) consisting of elements of zero norm we get a pre-Hilbert space, and its completion under the pre-inner product is denoted by $\clk$. We also define $\Gma_0 : E \raro \clb(h,\clk)$ by setting $$ (\Gma_0(\xi))w:=[\xi,w],$$ where $[\xi,w]$ represents the equivalence class of $(\xi,w)$ in $\cls \equiv \clv/\clv_0 \subseteq \clk$. That it is an isometry is verified by straightforward calculations. Next, we define $V_g$ on $\cls$ by $$ V_g [\xi,w]:=[\beta_g(\xi), u_g w],$$ and verify that it is indeed an isometry, and since its range clearly contains a total subset, $V_g$  extends to a unitary on $\clk$.
 Furthermore, $V_g
V_h=V_{gh}$ and $V_e={\rm Id}$ (where $e$ is the identity of $G$) on $\cls$,and hence
on the whole of $\clk.$ The strong continuity of $g \mapsto V_g$ is also easy
to see. Indeed, it is enough to prove that $g \mapsto V_g X$ is continuous for
any $X$ of the form $[\xi, v]$, $ \xi \in E, v\in h.$ But $\|V_g([\xi, v])-[\xi,
v]\|^2=2  \langle  [\xi, v], [\xi, v]  \rangle  -  \langle  V_g([\xi,
v]),[\xi, v]  \rangle  -  \langle  [\xi, v], V_g([\xi, v])  \rangle  ,$ and we have,
 $  \langle  V_g([\xi, v]),[\xi, v]  \rangle  -  \langle  [\xi, v], [\xi, v]  \rangle  =   \langle (u_g v-v)\langle \beta_g(\xi), \xi \rangle
v  \rangle  +  \langle v, \langle  (\beta_g(\xi)-\xi), \xi \rangle v  \rangle  .$   By assumption $ \lim_{g \raro e}
(\beta_g(\xi)- \xi)=0$   in the norm topology of $E$, so
 $  \langle  (\beta_g(\xi)-\xi)v,\xi v  \rangle   \raro 0$ as $g \raro e$. Furthermore,
  $g \mapsto u_g v$ is continuous. This completes the proof of strong continuity of $V_g$.
  \qed

 In view of the above result, we assume without loss of genetrality that $E \subset \clb(h,\clk)$ (with the natural Hilbert module structure inherited from that of $\clb(h,\clk)$), and $\beta_g(\cdot)=V_g \cdot u_g^{-1}$.
 Consider the strong operator closure of $\tilde{E}$ of $E$ in $\clb(h,\clk)$. It is a Hilbert von Neumann $\tilde{\cla}$ module (where $\tilde{\cla}$ is the weak closure of $\cla$ in $h$) . Moreover, the $G$-action $\beta_g=V_g \cdot u_g^{-1}$ can be extended to the whole of $\clb(h,\clk)$, and denoted again by $\beta_g$. Clearly, this action leaves $\tilde{E}$ invariant, hence $(\tilde{E},\beta)$ is a Hilbert von Neumann $G$-$\tilde{\cla}$ module. Let us recall that  by  ${\tilde{E}}^\infty$ we denote the  locally convex  space of elements $\xi$  in $\tilde{E}$ such that $g \mapsto \beta_g(\xi)$ is $C^\infty$ in the strong operator topology of $\tilde{E}$.
 \bthm
  There exist a Hilbert space $k_0$, a unitary representation $w_g$ of $G$ in $k_0$ and an isometry  $\Sigma$  from $\clk$ to $h \ot k_0$ such that\\
  (i) $\Sigma$ is equivariant in the sense that $\Sigma V_g=(u_g \ot w_g) \Sigma$ for all $g$;\\
  (ii) $\Sigma \xi \in \cla \ot k_0$ for all $\xi \in E$. \\
  \ethm
 {\it Proof :}\\
 The statement (i) is contained in the Theorem 4.3.5 of \cite{book} (page 99).
 For proving (ii), we note that $E^\infty$ (w.r.t. the action $\beta$) is mapped by $\Sigma$ into $(\tilde{\cla} \ot k_0)^\infty$ (w.r.t. the action $\gamma_g:={\rm ad}_{u_g} \ot w_g$), and moreover, for $\xi \in E$, $g \mapsto \gamma_g(\Sigma(\xi)) =\Sigma \beta_g(\xi)$ is norm-continuous since $g \mapsto \beta_g(\xi)$ is so and $\Sigma$ is isometry. Thus, (ii) follows from Lemma \ref{mainlemma}.
 \qed

It follows from the above theorem that $E$ can be equivariantly
embedded in the trivial $G-\cla$ module $(\cla \ot k_0, \al \ot w
)$. In particular, we have that \bthm If a compact Lie group $G$ has
an ergodic action on a $C^*$-algebra $\cla$, then every $G-\cla$
module $(E,\beta)$ is embeddable. \ethm \vspace{2mm}
 {\bf
Acknowledgement :} The author would like to thank the anonymous
referee who pointed out a crucial mistake in the earlier version,
which has led to substantial revision of the paper.

{\bf Stat-Math Unit, Indian Statistical Institute}\\
{\bf 203, B. T. Road, Kolkata 700108, India.}\\
{\bf E-mail : goswamid@isical.ac.in \\
\end{document}